\documentclass[11pt,a4paper]{article}

\usepackage[T1]{fontenc}
\usepackage[utf8]{inputenc}
\usepackage{lmodern}
\usepackage{microtype}
\usepackage{amsmath,amssymb,amsthm,mathtools}
\usepackage{bm}
\usepackage{booktabs}
\usepackage{array}
\usepackage{enumitem}
\usepackage{geometry}
\usepackage{xcolor}
\usepackage{tocloft}
\usepackage[colorlinks,final,backref=page,hyperindex]{hyperref}

\allowdisplaybreaks
\hypersetup{
  pdfauthor={Haiying Li, Yijia Zhao, and Yongjin Li},
  pdftitle={Matrix Amplification Profiles and Operator Models for James and von Neumann--Jordan Constants}
}

\newtheorem{theorem}{Theorem}[section]
\newtheorem{proposition}[theorem]{Proposition}
\newtheorem{lemma}[theorem]{Lemma}
\newtheorem{corollary}[theorem]{Corollary}
\theoremstyle{definition}
\newtheorem{definition}[theorem]{Definition}
\newtheorem{example}[theorem]{Example}
\theoremstyle{remark}
\newtheorem{remark}[theorem]{Remark}

\newcommand{\R}{\mathbb{R}}

\newcommand{\CNJ}{C_{\mathrm{NJ}}}
\newcommand{\Cprime}{C'_{\mathrm{NJ}}}
\newcommand{\CZ}{C_{\mathrm Z}}
\newcommand{\CBM}{C_{\mathrm{BM}}}
\newcommand{\MA}{\mathfrak M}
\newcommand{\JA}{\mathfrak J}
\newcommand{\Gam}{\Gamma_2}
\newcommand{\Lam}{\Lambda_2}
\newcommand{\cond}{\kappa_2}

\DeclareMathOperator{\spanop}{span}
\DeclareMathOperator{\diag}{diag}

\title{\textbf{Matrix Amplification Profiles and Operator Models\\
for James and von Neumann--Jordan Constants}}

\author{Haiying Li, Yijia Zhao\\[2mm]
\small  School of Mathematics and Statistics,\\
\small  Henan Normal University,\\
\small  Xinxiang 453007\\
\small \texttt{lihaiying@htu.edu.cn(Haiying Li);zhaoyijia20259@stu.htu.edu.cn(Yijia Zhao)}\\[2mm]
Yongjin Li\\[2mm]
\small Department of Mathematics, Sun Yat-sen University,\\
\small Guangzhou 510275, Guangdong, China\\
\small \texttt{stslyj@mail.sysu.edu.cn}}
\date{}

\begin{document}
\maketitle

\begin{abstract}
For a real Banach space $X$, we compare the norm of the operator induced on
$X\oplus_2X$ by a scalar $2\times2$ matrix with the matrix spectral norm.
This gives a common operator model for the von Neumann--Jordan and James
constants.  We obtain sharp matrix bounds, exact formulas on spaces
containing $\ell_1^2$ or $\ell_\infty^2$, and a global profile that
characterizes Hilbert spaces and uniform nonsquareness.  We prove duality,
Hilbertian-sum and Banach--Mazur stability, calculate the profile of $L_p$
and $\ell_p$, and derive relations with classical geometric constants and
moduli.  The profile also yields fixed point criteria for nonexpansive
mappings.
\end{abstract}




\section{Introduction}

Geometric constants compress a nonlinear feature of the unit ball into a
number that can be compared across Banach spaces.  Two of the most useful
examples are the James constant \cite{GaoLau1990}
\begin{equation}\label{eq:J-intro}
 J(X)=\sup_{x,y\in S_X}\min\{\|x+y\|,\|x-y\|\}
\end{equation}
and the von Neumann--Jordan constant
\begin{equation}\label{eq:CNJ-intro}
 \CNJ(X)=\sup_{(x,y)\ne(0,0)}
 \frac{\|x+y\|^2+\|x-y\|^2}
 {2(\|x\|^2+\|y\|^2)}.
\end{equation}
The first records whether two unit vectors can be simultaneously far apart
after addition and subtraction; the second measures the largest failure of
the parallelogram identity.  Their common algebraic core is the Hadamard
transformation
\[
 (x,y)\longmapsto(x+y,x-y).
\]
Nevertheless, the two constants are usually introduced and investigated
separately: $J(X)$ uses a minimum on the bisphere, whereas $\CNJ(X)$ uses a
quadratic quotient on the whole product space.

The aim of this paper is to place both constants in a transparent operator
framework.  Given $A\in M_2(\R)$, let $T_A^X$ act on the Hilbertian sum
$X\oplus_2X$ by applying the scalar matrix to the two $X$-valued
coordinates.  We compare its norm with the Euclidean spectral norm of $A$.
This comparison is trivial when $X$ is Hilbert, but it detects the geometry
of a general norm.  The Hadamard matrix recovers $\CNJ(X)$ exactly \cite{KatoTakahashi1997}, and its
two coordinate outputs recover $J(X)$ by a bottleneck aggregation.

There are two points of caution.  First, an orthogonal scalar matrix need
not induce an isometry of $X\oplus_2X$; the assertion is true for every
orthogonal matrix precisely in the Hilbertian situation.  Second, a fixed
non-scalar matrix cannot in general characterize inner-product spaces.
Indeed, every sign-balanced matrix has normalized amplification one on
\emph{every} Banach space.  These observations lead naturally to the global
profile $\Gam(X)$, obtained by taking the supremum over all scalar matrices,
and to a condition-number-restricted profile $\Lam(X;K)$.
For this definition, however, the two profiles coincide for every
$K\ge1$: the global supremum is already attained among orthogonal
matrices.

Our main contributions are as follows.
\begin{enumerate}[label=(\roman*),leftmargin=2.2em]
\item We establish $1\le \MA_A(X)\le
\|A^{\mathrm{abs}}\|_2/\|A\|_2\le\sqrt2$ and characterize the sign-balanced class on
which the lower bound is always attained.
\item If $X$ contains $\ell_1^2$ or $\ell_\infty^2$ isometrically, then the
upper bound is attained for every $A$.
\item We prove the two sharp geometric characterizations
\[
\Gam(X)=1\iff X\text{ is Hilbert},\qquad
\Gam(X)<\sqrt2\iff X\text{ is uniformly nonsquare}.
\]
\item We obtain the exact profiles
\[
 \Gam(\ell_p)=\Gam(L_p)=2^{|1/p-1/2|},\qquad 1\le p\le\infty,
\]
with two disjointly supported nonzero functions in $L_p$
and dimension at least two in the sequence-space case. We also
obtain an interpolation bound for each individual matrix.
\item We connect the operator constants with classical geometric constants,
duality, Hilbertian sums, the Banach--Mazur distance, moduli, and the fixed
point property for nonexpansive mappings.
\end{enumerate}

The paper is organized so that the operator construction is introduced
before any new profile.  The classical constants then appear as special
responses of one concrete matrix, making the comparison between their
usual proof techniques explicit.

\section{Preliminaries}

Throughout, $X$ is a nonzero real Banach space, $B_X$ and $S_X$ denote its
closed unit ball and unit sphere, respectively, and $\|A\|_2$ is the
Euclidean spectral norm of a real $2\times2$ matrix.  We equip
\[
 X\oplus_2X=X\times X
\quad\text{with}\quad
 \|(x,y)\|_{2,X}=(\|x\|^2+\|y\|^2)^{1/2}.
\]
We use $L_p$ for Lebesgue spaces and $\ell_p$ for sequence spaces
throughout; the notation is kept distinct also at the endpoints
$p=1,\infty$.

For a real matrix $A=(a_{ij})$, its \emph{entrywise absolute-value matrix} is
\[
 A^{\mathrm{abs}}=(|a_{ij}|)_{i,j=1}^{2}.
\]
Each entry of $A^{\mathrm{abs}}$ is the scalar absolute value of the
corresponding entry of $A$. We denote the spectral norms of these two
matrices by $\|A^{\mathrm{abs}}\|_2$ and $\|A\|_2$, respectively.
The positive semidefinite factor in the polar decomposition of $A$ is
written $(A^TA)^{1/2}$.

We use the order-two Hadamard matrix \cite[p.~101]{SaitoKomuroTanaka2017}
\[
 H=\begin{pmatrix}1&1\\1&-1\end{pmatrix}.
\]
Its action on a pair of vectors is $(x,y)\mapsto(x+y,x-y)$.
Since $H^TH=2I_2$, the normalized matrix $H/\sqrt2$ is orthogonal.

Besides \eqref{eq:J-intro} and \eqref{eq:CNJ-intro}, we use the modified
von Neumann--Jordan constant \cite[Eq.~(1)]{AlonsoMartinPapini2008}
\begin{equation}\label{eq:Cprime}
 \Cprime(X)=\sup_{x,y\in S_X}
 \frac{\|x+y\|^2+\|x-y\|^2}{4}
\end{equation}
and the Zb\u{a}ganu constant \cite{Zbaganu2002}
\begin{equation}\label{eq:CZ}
 \CZ(X)=\sup_{(x,y)\ne(0,0)}
 \frac{\|x+y\|\,\|x-y\|}{\|x\|^2+\|y\|^2}.
\end{equation}
By the arithmetic--geometric mean inequality,
\begin{equation}\label{eq:CZ-CNJ}
 \CZ(X)\le \CNJ(X).
\end{equation}
Also $1\le\Cprime(X)\le\CNJ(X)\le2$.  The equality
$\CNJ(X)=1$ is equivalent to the parallelogram identity and hence to the
existence of an inner product generating the norm
\cite{JordanVonNeumann1935,KatoTakahashi1997}.

Recall that $X$ is \emph{uniformly nonsquare} if there is $\eta>0$ such
that, for all $x,y\in B_X$,
\[
 \min\{\|x+y\|,\|x-y\|\}\le2-\eta.
\]
Equivalently,
\begin{equation}\label{eq:uns-equiv}
 J(X)<2\quad\Longleftrightarrow\quad \CNJ(X)<2.
\end{equation}
For $\dim X\ge2$, the supremum in the definition of $J(X)$ is unchanged
if $S_X$ is replaced by $B_X$ \cite[p.~280, property~(iv)]{KatoMaligrandaTakahashi2001}.
Thus $J(X)<2$ is equivalent to uniform nonsquareness, and the
equivalence with $\CNJ(X)<2$ is classical
\cite{James1964,TakahashiKato1998}. In dimension one, $J(X)=0$ and
$\CNJ(X)=1$. For $x=au$ and $y=bu$ in $B_X$, where $\|u\|=1$,
the smaller of $\|x+y\|$ and $\|x-y\|$ is $\bigl||a|-|b|\bigr|\le1$,
so uniform nonsquareness also holds in this case.

The modulus of convexity and modulus of smoothness are
\cite{Figiel1976}
\begin{align}
 \delta_X(\varepsilon)
 &=\inf\left\{1-\left\|\frac{x+y}{2}\right\|:
 x,y\in S_X,\ \|x-y\|\ge\varepsilon\right\},\label{eq:delta}\\
 \rho_X(t)
 &=\sup_{x,y\in S_X}
 \left(\frac{\|x+ty\|+\|x-ty\|}{2}-1\right).
 \label{eq:rho}
\end{align}
For nonzero $X$ we use the dimension-free form of the James--modulus
identity:
\begin{equation}\label{eq:J-delta}
 J(X)=\sup\left(\{0\}\cup\left\{\varepsilon\in(0,2):
 \delta_X(\varepsilon)\le1-\frac{\varepsilon}{2}\right\}\right).
\end{equation}
For $\dim X\ge2$ this is the classical formula
\cite[p.~280, Eq.~(2.2)]{KatoMaligrandaTakahashi2001}; the extra point $0$
also covers $\dim X=1$, when $J(X)=0$ and the positive part of the set
is empty.

For a two-dimensional normed space $E$, let
$d_{\mathrm{BM}}(E,\ell_2^2)$ be its Banach--Mazur distance from the Euclidean
plane.  Following Stern \cite[Definition~3.2]{Stern2015}, we set
\begin{equation}\label{eq:CBM}
 \CBM(X)=\sup\{d_{\mathrm{BM}}(E,\ell_2^2)^2:
 E\subset X,\ \dim E=2\},
\end{equation}
with value $1$ when $\dim X=1$.  John's theorem gives
$1\le\CBM(X)\le2$ \cite{John1948}.

\section{Matrix-induced operators}

\begin{definition}\label{def:MA}
Let $A=(a_{ij})\in M_2(\R)$.  Define
\[
 T_A^X:X\oplus_2X\longrightarrow X\oplus_2X,
 \qquad
 T_A^X(x,y)=
 (a_{11}x+a_{12}y,a_{21}x+a_{22}y).
\]
Its operator norm is given by
\[
 \|T_A^X\|
 =\sup_{\substack{x,y\in X\\(x,y)\ne(0,0)}}
 \frac{\|T_A^X(x,y)\|_{2,X}}{\|(x,y)\|_{2,X}}.
\]
For $A\ne0$, the \emph{normalized matrix amplification} of $X$ at $A$ is
\begin{equation}\label{eq:MA-def}
 \MA_A(X)=\frac{\|T_A^X\|}{\|A\|_2}.
\end{equation}
\end{definition}

The normalization removes scalar size: $\MA_{\lambda A}(X)=\MA_A(X)$ for
$\lambda\ne0$.  It also makes the Hilbert-space value equal to one.

\begin{theorem}\label{thm:universal}
For every nonzero real Banach space $X$ and every nonzero $A\in M_2(\R)$,
\begin{equation}\label{eq:universal}
 1\le\MA_A(X)\le
 \frac{\|A^{\mathrm{abs}}\|_2}{\|A\|_2}\le\sqrt2.
\end{equation}
All three bounds are sharp.
\end{theorem}

\begin{proof}
Fix $u\in S_X$ and consider the subspace
$V=\{(su,tu):s,t\in\R\}\subset X\oplus_2X$. For $s,t\in\R$,
\[
 T_A^X(su,tu)=((a_{11}s+a_{12}t)u,(a_{21}s+a_{22}t)u).
\]
Since $\|u\|=1$, the definition of the norm on $X\oplus_2X$ gives
\begin{align*}
 \|T_A^X(su,tu)\|_{2,X}
 &=\Bigl(|a_{11}s+a_{12}t|^2\|u\|^2
        +|a_{21}s+a_{22}t|^2\|u\|^2\Bigr)^{1/2}\\
 &=\left\|A\begin{pmatrix}s\\t\end{pmatrix}\right\|_2.
\end{align*}
Moreover,
\[
 \|(su,tu)\|_{2,X}
 =\bigl(s^2\|u\|^2+t^2\|u\|^2\bigr)^{1/2}
 =\left\|\begin{pmatrix}s\\t\end{pmatrix}\right\|_2.
\]
The operator norm is the supremum over all nonzero vectors in
$X\oplus_2X$. Restricting to the nonzero vectors $(su,tu)\in V$ yields
\[
 \|T_A^X\|
 \ge\sup_{(s,t)\ne(0,0)}
 \frac{\|T_A^X(su,tu)\|_{2,X}}{\|(su,tu)\|_{2,X}}
 =\sup_{(s,t)\ne(0,0)}
 \frac{\left\|A\begin{pmatrix}s\\t\end{pmatrix}\right\|_2}
      {\left\|\begin{pmatrix}s\\t\end{pmatrix}\right\|_2}
 =\|A\|_2.
\]

Next, for arbitrary $(x,y)\in X\oplus_2X$, put
$r=(\|x\|,\|y\|)^T$. The triangle inequality gives
\[
 \begin{pmatrix}\|a_{11}x+a_{12}y\|\\
 \|a_{21}x+a_{22}y\|\end{pmatrix}\le A^{\mathrm{abs}} r
\]
coordinatewise. All coordinates are nonnegative, so
\[
 \|T_A^X(x,y)\|_{2,X}
 \le\|A^{\mathrm{abs}} r\|_2\le\|A^{\mathrm{abs}}\|_2\|r\|_2
 =\|A^{\mathrm{abs}}\|_2\|(x,y)\|_{2,X}.
\]
Thus $\|T_A^X\|\le\|A^{\mathrm{abs}}\|_2$.

Finally, let $s_1\ge s_2\ge0$ be the singular values of $A$, and let
$\|B\|_{\mathrm F}=(\sum_{i,j}|b_{ij}|^2)^{1/2}$ denote the
Frobenius norm. Since $A\ne0$, $s_1=\|A\|_2>0$, and
\[
 \frac{\|A^{\mathrm{abs}}\|_2}{\|A\|_2}
 \le\frac{\|A^{\mathrm{abs}}\|_{\mathrm F}}{\|A\|_2}
 =\frac{\|A\|_{\mathrm F}}{\|A\|_2}
 =\sqrt{1+(s_2/s_1)^2}\le\sqrt2.
\]
These inequalities prove \eqref{eq:universal}.

The lower bound is attained on every nonzero Hilbert space, where
the singular-value decomposition gives $\|T_A^X\|=\|A\|_2$.
For the upper bounds, use the Hadamard matrix $H$ and take $X=\ell_1^2$.
Then $\|H\|_2=\sqrt2$ and $\|H^{\mathrm{abs}}\|_2=2$. For the standard unit vectors,
\[
 \frac{\|T_H^X(e_1,e_2)\|_{2,X}}{\|(e_1,e_2)\|_{2,X}}
 =\frac{\sqrt{2^2+2^2}}{\sqrt2}=2.
\]
Together with the upper estimate, this gives
$\MA_H(X)=\|H^{\mathrm{abs}}\|_2/\|H\|_2=\sqrt2$, so both upper bounds are sharp.
\end{proof}

\begin{remark}\label{rem:orthogonal-warning}
An orthogonal matrix $Q$ satisfies $\|Q\|_2=1$, but $T_Q^X$ need not
be an isometry. For
\[
 Q=\frac1{\sqrt2}\begin{pmatrix}1&1\\1&-1\end{pmatrix},
\]
the definition of $\CNJ(X)$ gives
\[
 \|T_Q^X\|^2=\sup_{(x,y)\ne(0,0)}
 \frac{\|x+y\|^2+\|x-y\|^2}{2(\|x\|^2+\|y\|^2)}
 =\CNJ(X).
\]
For $X=\ell_1^2$, the preceding proof gives $\|T_Q^X\|=\sqrt2>1$.
Thus scalar orthogonality does not imply isometry on $X\oplus_2X$.
Every orthogonal matrix induces an isometry when $X$ is Hilbert;
conversely, isometry for this particular normalized Hadamard matrix
already forces $\CNJ(X)=1$ and hence the Hilbert-space property.
\end{remark}

It is convenient to name the purely matrix-theoretic upper bound.
\begin{definition}
For $A\ne0$, set
\begin{equation}\label{eq:frustration}
 m(A)=\frac{\|A^{\mathrm{abs}}\|_2}{\|A\|_2}.
\end{equation}
We call $m(A)$ the \emph{sign-frustration ratio} of $A$.
Then $1\le\MA_A(X)\le m(A)\le\sqrt2$.
\end{definition}

The ratio $m(A)$ depends only on the matrix $A$ and gives
a universal upper bound for $\MA_A(X)$. If $A$ is sign-balanced,
Proposition~\ref{prop:balanced} gives $\MA_A(X)=m(A)=1$.
For the Hadamard matrix $H$, $m(H)=\sqrt2$; this upper bound is
attained on the spaces in Theorem~\ref{thm:extreme-exact}, but need not
be attained on every space.

\section{Sign-balanced matrices and exact endpoint spaces}

Multiplying a row of $A$ by $-1$ reverses the sign of one output
coordinate of $T_A^X$. Multiplying a column by $-1$ reverses the sign of
one input coordinate before $T_A^X$ is applied. Both operations preserve
$\|T_A^X\|$ and $\|A\|_2$, and hence preserve $\MA_A(X)$.
This lets us characterize the nonzero matrices for which the lower
bound $\MA_A(X)=1$ holds on every nonzero real Banach space. We also
show that if $X$ contains an isometric copy of $\ell_1^2$ or
$\ell_\infty^2$, then $\MA_A(X)=m(A)$ for every nonzero $A$.

\begin{definition}
A matrix $A\in M_2(\R)$ is \emph{sign-balanced} if there exist diagonal
sign matrices $D_1,D_2$, with diagonal entries in $\{-1,1\}$, such that
$D_1AD_2$ is entrywise nonnegative.
\end{definition}

For a $2\times2$ matrix with no zero entries, sign balance is equivalent to
$a_{11}a_{12}a_{21}a_{22}>0$. Indeed, the product is invariant under
row and column sign changes. By choosing the signs of the first row,
the second column and the second row in that order, we can make
$a_{11}$, $a_{12}$ and $a_{21}$ nonnegative. When all entries are
nonzero, the sign of the remaining entry is then the sign of the product.
If $A$ has a zero entry, label that entry $a_{22}$ by reindexing the
rows and columns. The same choice of signs makes every nonzero entry
positive, so $A$ is sign-balanced.
More explicitly, for diagonal sign matrices $D_1,D_2$ one has
\begin{equation}\label{eq:sign-orbit}
 T_{D_1AD_2}^X=T_{D_1}^X T_A^X T_{D_2}^X,
 \qquad
 \|D_1AD_2\|_2=\|A\|_2.
\end{equation}
Both $T_{D_1}^X$ and $T_{D_2}^X$ are isometries.  Consequently, normalized
amplification is constant on every row--column sign orbit, and the study of
sign-balanced matrices is exactly equivalent to the study of their
entrywise nonnegative representatives.  The same observation applies to
row and column permutations.

\begin{proposition}\label{prop:balanced}
If $A\ne0$ is sign-balanced, then
\begin{equation}\label{eq:balanced}
 \MA_A(X)=1
\end{equation}
for every real Banach space $X$.
\end{proposition}

\begin{proof}
Since $A$ is sign-balanced, there exist diagonal sign matrices $D_1,D_2$,
with diagonal entries in $\{-1,1\}$, such that $B=D_1AD_2$ is
entrywise nonnegative. By \eqref{eq:sign-orbit},
\[
 T_B^X=T_{D_1}^XT_A^XT_{D_2}^X.
\]
The operators $T_{D_1}^X$ and $T_{D_2}^X$ are surjective isometries on
$X\oplus_2X$, so $\|T_A^X\|=\|T_B^X\|$. Orthogonal invariance of the
spectral norm gives $\|B\|_2=\|A\|_2$; hence
\[
 \MA_A(X)=\MA_B(X).
\]
For the entrywise nonnegative matrix $B$, Theorem~\ref{thm:universal}
yields
\[
 \|B\|_2\le\|T_B^X\|\le\|B^{\mathrm{abs}}\|_2=\|B\|_2.
\]
Thus $\|T_B^X\|=\|B\|_2$. Therefore $\MA_B(X)=1$, which implies
$\MA_A(X)=1$.
\end{proof}

\begin{remark}
Proposition~\ref{prop:balanced} shows why a statement of the form
``$\MA_A(X)=1$ for one fixed non-scalar matrix $A$ if and only if $X$ is
Hilbert'' is false without a strong hypothesis on $A$.  Nonnegative
matrices, triangular matrices with one zero entry, and all their row and
column sign changes have value one on every Banach space.
\end{remark}

The opposite endpoint is also completely computable.

\begin{theorem}
\label{thm:extreme-exact}
Suppose that $X$ contains an isometric copy of $\ell_1^2$ or of
$\ell_\infty^2$.  Then, for every nonzero $A\in M_2(\R)$,
\begin{equation}\label{eq:extreme-exact}
 \MA_A(X)=m(A)=\frac{\|A^{\mathrm{abs}}\|_2}{\|A\|_2}.
\end{equation}
\end{theorem}

\begin{proof}
By Theorem~\ref{thm:universal},
$\MA_A(X)\le\|A^{\mathrm{abs}}\|_2/\|A\|_2=m(A)$.
It remains to prove $\|T_A^X\|\ge\|A^{\mathrm{abs}}\|_2$.
Set $G=\bigl(A^{\mathrm{abs}}\bigr)^TA^{\mathrm{abs}}$. For a nonzero
$v\in\R^2$, the Rayleigh quotient of $G$ at $v$ is $v^TGv/(v^Tv)$.
Its maximum is $\lambda_{\max}(G)=\|A^{\mathrm{abs}}\|_2^2$, and any
unit vector attaining this maximum is an eigenvector for that eigenvalue.
Since $G$ has nonnegative entries, taking the absolute value of each
coordinate of a maximizing unit vector preserves its norm and cannot
decrease the quotient. We can therefore choose a nonnegative maximizing
unit vector $r=(r_1,r_2)^T$, for which
\[
 \|A^{\mathrm{abs}}r\|_2^2=r^TGr=\|A^{\mathrm{abs}}\|_2^2.
\]

\emph{Case 1: $X$ contains an isometric copy of $\ell_1^2$.}
Take $e_1,e_2\in S_X$ corresponding to the standard basis of this copy,
so $\|se_1+te_2\|=|s|+|t|$ for all real $s,t$, and set
\[
 x=r_1e_1,\qquad y=r_2e_2.
\]
Then $\|(x,y)\|_{2,X}=(r_1^2+r_2^2)^{1/2}=1$, and, for $i=1,2$,
\[
 \|a_{i1}x+a_{i2}y\|
 =\|a_{i1}r_1e_1+a_{i2}r_2e_2\|
 =|a_{i1}|r_1+|a_{i2}|r_2.
\]
Consequently,
\[
 \|T_A^X(x,y)\|_{2,X}
 =\left((|a_{11}|r_1+|a_{12}|r_2)^2
       +(|a_{21}|r_1+|a_{22}|r_2)^2\right)^{1/2}
 =\|A^{\mathrm{abs}} r\|_2=\|A^{\mathrm{abs}}\|_2.
\]
Hence $\|T_A^X\|\ge\|A^{\mathrm{abs}}\|_2$.

\emph{Case 2: $X$ contains an isometric copy of $\ell_\infty^2$.}
Take $e_1,e_2$ corresponding to the standard basis and use coordinates
in this copy, with its inherited norm denoted by $\|\cdot\|_\infty$.
Set
\[
 x=r_1(\operatorname{sgn}a_{11},\operatorname{sgn}a_{21}),\qquad
 y=r_2(\operatorname{sgn}a_{12},\operatorname{sgn}a_{22}).
\]
In this construction, the sign of a zero entry is chosen arbitrarily
from $\{-1,1\}$. Then $\|x\|_\infty=r_1$, $\|y\|_\infty=r_2$,
and $\|(x,y)\|_{2,X}=1$. The first coordinate of $a_{11}x+a_{12}y$
and the second coordinate of $a_{21}x+a_{22}y$ give
\[
 \|a_{11}x+a_{12}y\|_\infty\ge|a_{11}|r_1+|a_{12}|r_2,
 \qquad
 \|a_{21}x+a_{22}y\|_\infty\ge|a_{21}|r_1+|a_{22}|r_2.
\]
The triangle inequality gives the reverse estimates, so both are
equalities. Therefore
\[
 \|T_A^X(x,y)\|_{2,X}=\|A^{\mathrm{abs}} r\|_2=\|A^{\mathrm{abs}}\|_2,
\]
and again $\|T_A^X\|\ge\|A^{\mathrm{abs}}\|_2$.

Combining this lower bound in both cases with Theorem~\ref{thm:universal}
gives $\|T_A^X\|=\|A^{\mathrm{abs}}\|_2$. Dividing by $\|A\|_2>0$ yields
\[
 \MA_A(X)=\frac{\|T_A^X\|}{\|A\|_2}
 =\frac{\|A^{\mathrm{abs}}\|_2}{\|A\|_2}=m(A).
\]
\end{proof}

Consequently, $\MA_A(X)=1$ for every nonzero real Banach space $X$ if and
only if $A\ne0$ is sign-balanced. Indeed, the sign-balanced case is
Proposition~\ref{prop:balanced}. Otherwise all four entries of $A$ are
nonzero and $a_{11}a_{12}a_{21}a_{22}<0$. The matrices $A^TA$ and
$\bigl(A^{\mathrm{abs}}\bigr)^TA^{\mathrm{abs}}$ have equal diagonal entries, while the absolute value of
the off-diagonal entry of the latter is strictly larger. The eigenvalue
formula for a symmetric $2\times2$ matrix therefore gives
$\|A^{\mathrm{abs}}\|_2>\|A\|_2$. Theorem~\ref{thm:extreme-exact} yields
$\MA_A(\ell_1^2)=\|A^{\mathrm{abs}}\|_2/\|A\|_2>1$.

\begin{corollary}\label{cor:classical-endpoints}
Formula \eqref{eq:extreme-exact} holds for $\ell_1$, $\ell_\infty$, $c_0$,
and for $L_1(\mu)$ or $L_\infty(\mu)$ whenever the measure space contains
two disjoint sets of positive finite measure.  It also holds for $C(K)$
when the compact Hausdorff space $K$ has at least two points.
\end{corollary}

\begin{proof}
The sequence spaces contain the required two-dimensional coordinate
subspace.  Disjoint normalized characteristic functions give $\ell_1^2$
inside $L_1$ and $\ell_\infty^2$ inside $L_\infty$.  For two distinct
points of $K$, normality and Urysohn's lemma give two norm-one continuous
functions with disjoint supports, whose span is isometric to
$\ell_\infty^2$.
\end{proof}

\begin{example}\label{ex:At}
The upper bound $m(A)$ can vary even among matrices that become
orthogonal after normalization. For $t\ge0$, let
\[
 A_t=\begin{pmatrix}1&t\\ t&-1\end{pmatrix}.
\]
Since $A_t^TA_t=(1+t^2)I_2$, we have $\|A_t\|_2=\sqrt{1+t^2}$ and
$A_t/\sqrt{1+t^2}$ is orthogonal. Also, $\|A_t^{\mathrm{abs}}\|_2=1+t$.
Thus
\begin{equation}\label{eq:At}
 1\le\MA_{A_t}(X)\le m(A_t)=\frac{1+t}{\sqrt{1+t^2}}.
\end{equation}

At $t=0$, the map is $(x,y)\mapsto(x,-y)$, so $\MA_{A_0}(X)=1$
on every nonzero real Banach space. At $t=1$, we have $A_1=H$ and
$m(H)=\sqrt2$. If $X$ contains an isometric copy of $\ell_1^2$ or
$\ell_\infty^2$, Theorem~\ref{thm:extreme-exact} gives
$\MA_{A_t}(X)=m(A_t)$ for every $t\ge0$. On each such space,
varying $t$ from $0$ to $1$ therefore gives every amplification value
between $1$ and $\sqrt2$, including both endpoints.
\end{example}

\section{Global and condition-restricted profiles}

\begin{definition}\label{def:profiles}
The \emph{global $2\times2$ matrix amplification profile} is
\begin{equation}\label{eq:Gamma}
 \Gam(X)=\sup_{A\ne0}\MA_A(X).
\end{equation}
For $K\ge1$, define the condition-restricted profile
\begin{equation}\label{eq:Lambda}
 \Lam(X;K)=\sup\{\MA_A(X): A\text{ is invertible and }
 \cond(A)\le K\},
\end{equation}
where $\cond(A)=\|A\|_2\|A^{-1}\|_2$.
\end{definition}

The bounds $1\le\Lam(X;K)\le\Gam(X)\le\sqrt2$ follow from
Theorem~\ref{thm:universal} and the definitions. The two profiles are
equal because $\Gam(X)$ is attained by an orthogonal matrix, whose
condition number is $1$.

\begin{lemma}\label{lem:continuity}
For all $A,B\in M_2(\R)$,
\begin{equation}\label{eq:operator-Lipschitz}
 \|T_A^X-T_B^X\|\le\sqrt2\,\|A-B\|_2.
\end{equation}
Consequently, $A\mapsto\|T_A^X\|$ is $\sqrt2$-Lipschitz,
$A\mapsto\MA_A(X)$ is continuous on $M_2(\R)\setminus\{0\}$, and
\begin{equation}\label{eq:Gamma-maximum}
 \Gam(X)=\max_{\|A\|_2=1}\|T_A^X\|.
\end{equation}
Moreover, let
$O(2)=\{Q\in M_2(\R):Q^TQ=I_2\}$ denote the real orthogonal group.
For every $K\ge1$,
\begin{equation}\label{eq:Lambda-limit}
 \Lam(X;K)=\Gam(X)=\max_{Q\in O(2)}\|T_Q^X\|.
\end{equation}
\end{lemma}

\begin{proof}
By the definition of $T_A^X$, the operator depends linearly on the
matrix, and therefore $T_A^X-T_B^X=T_{A-B}^X$. By
Theorem~\ref{thm:universal},
\[
 \|T_A^X-T_B^X\|=\|T_{A-B}^X\|
 \le\|(A-B)^{\mathrm{abs}}\|_2
 \le\sqrt2\,\|A-B\|_2.
\]
When $A=B$, the inequalities are immediate; otherwise the theorem
applies to the nonzero matrix $A-B$. This proves
\eqref{eq:operator-Lipschitz}.

The reverse triangle inequality for the operator norm, together with
\eqref{eq:operator-Lipschitz}, gives
\[
 \bigl|\|T_A^X\|-\|T_B^X\|\bigr|
 \le\|T_A^X-T_B^X\|
 \le\sqrt2\,\|A-B\|_2.
\]
Hence $A\mapsto\|T_A^X\|$ is $\sqrt2$-Lipschitz and therefore
continuous. On $M_2(\R)\setminus\{0\}$,
\[
 \MA_A(X)=\frac{\|T_A^X\|}{\|A\|_2},\qquad \|A\|_2>0.
\]
As the denominator is a positive continuous function,
$A\mapsto\MA_A(X)$ is continuous on $M_2(\R)\setminus\{0\}$.

By definition, $\Gam(X)=\sup_{A\ne0}\MA_A(X)$. The homogeneity
$\MA_{\lambda A}(X)=\MA_A(X)$ for $\lambda\ne0$ allows us to
restrict this supremum to
\[
 \mathcal S_2=\{A\in M_2(\R):\|A\|_2=1\}.
\]
The set $\mathcal S_2$ is closed and bounded in the finite-dimensional
space $M_2(\R)$, hence compact. Therefore the continuous function
$A\mapsto\|T_A^X\|$ attains its maximum on $\mathcal S_2$, and
\[
 \Gam(X)=\sup_{A\in\mathcal S_2}\|T_A^X\|
 =\max_{A\in\mathcal S_2}\|T_A^X\|.
\]
This proves \eqref{eq:Gamma-maximum}.

To see why an orthogonal matrix attains the maximum, fix $A\in\mathcal S_2$. Since $\|A\|_2=1$,
write its singular value decomposition as $A=U\diag(1,t)V^T$, where $U,V\in O(2)$
and $0\le t\le1$. Let
\[
 Q_+=UV^T,\qquad Q_-=U\diag(1,-1)V^T.
\]
Both matrices are orthogonal and
\[
 A=\frac{1+t}{2}Q_++\frac{1-t}{2}Q_-.
\]
By linearity and convexity of the operator norm,
\[
 \|T_A^X\|\le\frac{1+t}{2}\|T_{Q_+}^X\|
 +\frac{1-t}{2}\|T_{Q_-}^X\|
 \le\max_{Q\in O(2)}\|T_Q^X\|.
\]
The preceding inequality holds for every $A\in\mathcal S_2$. Since
$O(2)$ is compact and $Q\mapsto\|T_Q^X\|$ is continuous, the maximum
over $O(2)$ exists. Combining the inequality with
\eqref{eq:Gamma-maximum}, we obtain
\[
 \Gam(X)=\max_{A\in\mathcal S_2}\|T_A^X\|
 \le\max_{Q\in O(2)}\|T_Q^X\|.
\]
On the other hand, $O(2)\subset\mathcal S_2$, and therefore
\[
 \max_{Q\in O(2)}\|T_Q^X\|
 \le\max_{A\in\mathcal S_2}\|T_A^X\|=\Gam(X).
\]
Consequently,
\[
 \Gam(X)=\max_{Q\in O(2)}\|T_Q^X\|.
\]
For every $Q\in O(2)$, its condition number satisfies
$\cond(Q)=\|Q\|_2\|Q^{-1}\|_2=1$. Hence every orthogonal matrix is
admissible in the definition of $\Lam(X;K)$ for $K\ge1$. Since
$\|Q\|_2=1$, we obtain
\[
 \Lam(X;K)\ge\sup_{Q\in O(2)}\MA_Q(X)
 =\max_{Q\in O(2)}\|T_Q^X\|=\Gam(X).
\]
Conversely, $\Lam(X;K)\le\Gam(X)$ because the matrices in the
definition of $\Lam(X;K)$ form a subset of the nonzero matrices.
Thus $\Lam(X;K)=\Gam(X)$, proving \eqref{eq:Lambda-limit}.
\end{proof}

Recall that
\begin{equation}\label{eq:H}
 H=\begin{pmatrix}1&1\\1&-1\end{pmatrix}.
\end{equation}

The identity $\MA_H(X)^2=\CNJ(X)$ is the normalized form of
Proposition~B(i) of Kato and Takahashi \cite[p.~1056]{KatoTakahashi1997}.
We recall the proof to compare $\MA_H(X)$ with $\Gam(X)$ and $\Lam(X;K)$.

\begin{proposition}\label{prop:Hadamard}
For every nonzero real Banach space $X$,
\begin{equation}\label{eq:Hadamard-CNJ}
 \MA_H(X)^2=\CNJ(X).
\end{equation}
Consequently,
\begin{equation}\label{eq:profile-chain-basic}
 \sqrt{\CNJ(X)}\le\Lam(X;K)=\Gam(X)\le\sqrt2
 \qquad(K\ge1).
\end{equation}
\end{proposition}

\begin{proof}
By Definition~\ref{def:MA}, the matrix $H$ induces the operator
$T_H^X(x,y)=(x+y,x-y)$. The norm on $X\oplus_2X$ therefore gives
\[
 \|T_H^X(x,y)\|_{2,X}^2=\|x+y\|^2+\|x-y\|^2,
 \qquad
 \|(x,y)\|_{2,X}^2=\|x\|^2+\|y\|^2.
\]
The quotients defining the operator norm are nonnegative, so squaring
their supremum is equivalent to taking the supremum of their squares.
Consequently,
\[
 \|T_H^X\|^2
 =\sup_{\substack{x,y\in X\\(x,y)\ne(0,0)}}
   \frac{\|x+y\|^2+\|x-y\|^2}{\|x\|^2+\|y\|^2}
 =2\CNJ(X),
\]
where the last equality follows from \eqref{eq:CNJ-intro}.
Since
\[
 H^TH=
 \begin{pmatrix}1&1\\1&-1\end{pmatrix}^{T}
 \begin{pmatrix}1&1\\1&-1\end{pmatrix}
 =2I_2,
\]
the two singular values of $H$ are $\sqrt2$, and hence
$\|H\|_2^2=2$. The normalization in \eqref{eq:MA-def} now yields
\[
 \MA_H(X)^2=\frac{\|T_H^X\|^2}{\|H\|_2^2}
 =\frac{2\CNJ(X)}{2}=\CNJ(X).
\]

The identity $H^TH=2I_2$ also gives $H^{-1}=\tfrac12 H^T$, so
$\|H^{-1}\|_2=1/\sqrt2$ and $\cond(H)=1$. Therefore $H$ is
admissible in \eqref{eq:Lambda} for every $K\ge1$. Since
$\MA_H(X)$ is nonnegative, \eqref{eq:Hadamard-CNJ} gives
\[
 \sqrt{\CNJ(X)}=\MA_H(X)\le\Lam(X;K).
\]
Lemma~\ref{lem:continuity} gives $\Lam(X;K)=\Gam(X)$.
Theorem~\ref{thm:universal} gives $\MA_A(X)\le\sqrt2$ for every
nonzero $A$. Taking the supremum over these matrices gives
$\Gam(X)\le\sqrt2$, which completes \eqref{eq:profile-chain-basic}.
\end{proof}

The lower bound $\Gam(X)=1$ characterizes Hilbert spaces, and a
strict gap below $\sqrt2$ characterizes uniform nonsquareness.

\begin{theorem}\label{thm:global-profile}
For every nonzero real Banach space $X$,
\begin{align}
 \Gam(X)=1
 &\quad\Longleftrightarrow\quad X\text{ is a Hilbert space},
 \label{eq:Gamma-Hilbert}\\
 \Gam(X)<\sqrt2
 &\quad\Longleftrightarrow\quad X\text{ is uniformly nonsquare}.
 \label{eq:Gamma-UNS}
\end{align}
Moreover, for every $K\ge1$,
\begin{equation}\label{eq:Lambda-UNS}
 \Lam(X;K)<\sqrt2
 \quad\Longleftrightarrow\quad X\text{ is uniformly nonsquare}.
\end{equation}
\end{theorem}

\begin{proof}
\leavevmode\par\nopagebreak[4]
\noindent\textbf{Hilbert space characterization \eqref{eq:Gamma-Hilbert}.}
\par\nopagebreak[4]
Suppose that $X$ is a Hilbert space. Let $Q=(q_{ij})\in O(2)$
and $x_1,x_2\in X$. Expanding the squared norms in the real inner
product and using $Q^TQ=I_2$, we obtain
\begin{align*}
 \|T_Q^X(x_1,x_2)\|_{2,X}^2
 &=\sum_{i=1}^2\Bigl\|\sum_{j=1}^2 q_{ij}x_j\Bigr\|^2\\
 &=\sum_{j,k=1}^2 (Q^TQ)_{jk}\langle x_j,x_k\rangle\\
 &=\|x_1\|^2+\|x_2\|^2.
\end{align*}
Thus $\|T_Q^X\|=1$ for every $Q\in O(2)$. The maximum formula
\eqref{eq:Lambda-limit} gives $\Gam(X)=1$.

Conversely, if $\Gam(X)=1$, then Proposition~\ref{prop:Hadamard}
gives
\[
 1\le\sqrt{\CNJ(X)}=\MA_H(X)\le\Gam(X)=1.
\]
Hence $\CNJ(X)=1$. By the characterization recalled in Section~2,
the norm satisfies the parallelogram identity and is induced by
an inner product. Completeness makes $X$ a Hilbert space.

\smallskip
\noindent\textbf{Uniform nonsquareness \eqref{eq:Gamma-UNS}.}
\par\nopagebreak[4]
By \eqref{eq:uns-equiv} and \eqref{eq:profile-chain-basic}, it remains
to show that $\Gam(X)=\sqrt2$ implies $\CNJ(X)=2$.

Suppose that $\Gam(X)=\sqrt2$.
By the maximum formula \eqref{eq:Gamma-maximum} in
Lemma~\ref{lem:continuity}, there is a matrix $A$ such that
$\|A\|_2=1$ and $\|T_A^X\|=\sqrt2$. The estimates in the proof of
Theorem~\ref{thm:universal} give
\[
 \sqrt2=\|T_A^X\|\le\|A^{\mathrm{abs}}\|_2
 \le\|A^{\mathrm{abs}}\|_{\mathrm F}=\|A\|_{\mathrm F}
 \le\sqrt2\,\|A\|_2=\sqrt2.
\]
Thus equality holds throughout the chain.

First, let $s_1\ge s_2\ge0$ be the singular values of $A$. Since
$s_1=\|A\|_2=1$ and $s_1^2+s_2^2=\|A\|_{\mathrm F}^2=2$,
we have $s_2=1$. Hence both eigenvalues of the symmetric matrix
$A^TA$ are one, so $A^TA=I_2$ and $A$ is orthogonal.

Next, let $t_1\ge t_2\ge0$ be the singular values of $A^{\mathrm{abs}}$.
The equality $\|A^{\mathrm{abs}}\|_2=\|A^{\mathrm{abs}}\|_{\mathrm F}$ gives
$t_1^2=t_1^2+t_2^2$. Thus $t_2=0$, and $t_1=\sqrt2>0$
shows that $A^{\mathrm{abs}}$ has rank one. Write
\[
 A^{\mathrm{abs}}=\begin{pmatrix}\alpha&\beta\\\gamma&\delta\end{pmatrix},
 \qquad \alpha,\beta,\gamma,\delta\ge0.
\]
Since $A$ is orthogonal, the rows and columns of $A^{\mathrm{abs}}$ have
Euclidean norm one. Therefore
\[
 \alpha^2+\beta^2=\gamma^2+\delta^2
 =\alpha^2+\gamma^2=\beta^2+\delta^2=1.
\]
Comparing the first sum with the third and fourth gives
$\beta=\gamma$ and $\alpha=\delta$, respectively. As $A^{\mathrm{abs}}$ has
rank one,
\[
 0=\det(A^{\mathrm{abs}})=\alpha\delta-\beta\gamma=\alpha^2-\beta^2.
\]
Thus $\alpha=\beta$, and $\alpha^2+\beta^2=1$ gives
$\alpha=\beta=\gamma=\delta=2^{-1/2}$.

Multiply the columns of $A$ by signs to make its first row
positive. Then choose the sign of the second row so that its
first entry is positive. These operations give diagonal sign
matrices $D_1,D_2$ such that
\[
 D_1AD_2=2^{-1/2}\begin{pmatrix}1&1\\1&\varepsilon\end{pmatrix},
 \qquad \varepsilon\in\{-1,1\}.
\]
The transformed rows are orthogonal, so $(1+\varepsilon)/2=0$.
It follows that $D_1AD_2=2^{-1/2}H$. By \eqref{eq:sign-orbit},
$T_{D_1AD_2}^X=T_{D_1}^XT_A^XT_{D_2}^X$.
Since both $T_{D_1}^X$ and $T_{D_2}^X$ are surjective isometries,
they preserve the operator norm under composition. Hence
\[
 \sqrt2=\|T_A^X\|=\|T_{2^{-1/2}H}^X\|
 =\frac{\|T_H^X\|}{\|H\|_2}=\MA_H(X).
\]
Proposition~\ref{prop:Hadamard} now gives $\CNJ(X)=\MA_H(X)^2=2$.
This proves the endpoint equivalence and hence \eqref{eq:Gamma-UNS}.

\smallskip
\noindent\textbf{Condition number restriction \eqref{eq:Lambda-UNS}.}
\par\nopagebreak[4]
For every $K\ge1$, Lemma~\ref{lem:continuity} gives
$\Lam(X;K)=\Gam(X)$. The conclusion follows from
\eqref{eq:Gamma-UNS}.
\end{proof}

\begin{remark}
Compactness gives a matrix that attains $\Gam(X)$. The equality case
$\Gam(X)=\sqrt2$ thus provides a specific matrix to which the equality
conditions in the norm estimates apply. Pointwise strict inequalities
$\MA_A(X)<\sqrt2$ alone would still allow the values to approach
$\sqrt2$ as $A$ varies.
\end{remark}

\begin{corollary}\label{cor:all-orthogonal}
The following are equivalent:
\begin{enumerate}[label=(\roman*),leftmargin=2em]
\item $X$ is a Hilbert space;
\item $T_Q^X$ is an isometry of $X\oplus_2X$ for every orthogonal
$Q\in M_2(\R)$;
\item $T_{2^{-1/2}H}^X$ is an isometry.
\end{enumerate}
\end{corollary}

\begin{proof}
$(i)\Rightarrow(ii)$.
Suppose that $X$ is a Hilbert space, and let
$Q=(q_{ij})\in M_2(\mathbb R)$ be orthogonal.
For every $(x,y)\in X\oplus_2 X$, we have
\begin{align*}
 \|T_Q^X(x,y)\|_{2,X}^2
 &= \|q_{11}x+q_{12}y\|^2
    +\|q_{21}x+q_{22}y\|^2\\
 &= (q_{11}^2+q_{21}^2)\|x\|^2
    +(q_{12}^2+q_{22}^2)\|y\|^2\\
 &\quad
    +2(q_{11}q_{12}+q_{21}q_{22})\langle x,y\rangle.
\end{align*}
Since $Q^TQ=I_2$,
\[
 q_{11}^2+q_{21}^2=1,\qquad
 q_{12}^2+q_{22}^2=1,\qquad
 q_{11}q_{12}+q_{21}q_{22}=0.
\]
Consequently,
\[
 \|T_Q^X(x,y)\|_{2,X}^2
 =\|x\|^2+\|y\|^2
 =\|(x,y)\|_{2,X}^2.
\]
Thus $T_Q^X$ is an isometry.

\medskip\noindent
$(ii)\Rightarrow(iii)$.
Since $H^TH=2I_2$, the matrix $Q=2^{-1/2}H$ is orthogonal.
Therefore condition $(ii)$ implies that $T_{2^{-1/2}H}^X$ is an isometry.

\medskip\noindent
$(iii)\Rightarrow(i)$.
Suppose that $T_{2^{-1/2}H}^X$ is an isometry. Then
$\|T_{2^{-1/2}H}^X\|=1$.
By linearity of the construction with respect to the matrix,
\[
 T_{2^{-1/2}H}^X=2^{-1/2}T_H^X.
\]
Hence $\|T_H^X\|=\sqrt2$. Since $\|H\|_2=\sqrt2$, it follows that
\[
 \MA_H(X)=\frac{\|T_H^X\|}{\|H\|_2}=1.
\]
By Proposition~\ref{prop:Hadamard},
\[
 \CNJ(X)=\MA_H(X)^2=1.
\]
The equality $\CNJ(X)=1$ implies that the norm satisfies the
parallelogram identity. By the Jordan--von Neumann theorem, the norm is
induced by an inner product. Since $X$ is complete, $X$ is a Hilbert space.
\end{proof}

\section{An operator formulation of the James constant}

The same operator $T_H^X(x,y)=(x+y,x-y)$ represents both classical
constants. Proposition~\ref{prop:Hadamard} expresses $\CNJ(X)$ in terms
of its operator norm, which uses the sum of the squared output norms.
For $J(X)$, the input vectors are both unit vectors, and the quantity
to maximize is the smaller of the two output norms. The following
definition applies that rule to any matrix $A$.

\begin{definition}\label{def:JA}
For $A=(a_{ij})\in M_2(\R)$, define
\begin{equation}\label{eq:JA}
 \JA_A(X)=\sup_{x,y\in S_X}
 \min\{\|a_{11}x+a_{12}y\|,
          \|a_{21}x+a_{22}y\|\}.
\end{equation}
When $A\ne0$, put
$\widehat{\JA}_A(X)=\JA_A(X)/\|A\|_2$.
\end{definition}

For $A=H$, the two outputs are $x+y$ and $x-y$, so the
definition of $\JA_H(X)$ is exactly the definition of $J(X)$.
The next proposition compares this minimum of the two output norms
with the operator norm.

\begin{proposition}\label{prop:energy-bottle}
For every $A\ne0$,
\begin{equation}\label{eq:energy-bottle}
 \widehat{\JA}_A(X)\le\MA_A(X)\le\Gam(X).
\end{equation}
For the Hadamard matrix,
\begin{equation}\label{eq:J-H}
 \JA_H(X)=J(X),\qquad
 \widehat{\JA}_H(X)=\frac{J(X)}{\sqrt2}.
\end{equation}
Consequently,
\begin{equation}\label{eq:J-CNJ-Gamma}
 \frac{J(X)}{\sqrt2}
 \le\sqrt{\CNJ(X)}
 \le\Gam(X)\le\sqrt2.
\end{equation}
\end{proposition}

\begin{proof}
For $x,y\in S_X$, write $T_A^X(x,y)=(u_1,u_2)$. The smaller of
$\|u_1\|^2$ and $\|u_2\|^2$ is at most their arithmetic mean.
Taking square roots and using $\|(x,y)\|_{2,X}=\sqrt2$, we obtain
\[
 \min\{\|u_1\|,\|u_2\|\}
 \le\left(\frac{\|u_1\|^2+\|u_2\|^2}{2}\right)^{1/2}
 =\frac{\|T_A^X(x,y)\|_{2,X}}{\|(x,y)\|_{2,X}}
 \le\|T_A^X\|.
\]
Taking the supremum over $x,y\in S_X$ and dividing by $\|A\|_2$
gives $\widehat{\JA}_A(X)\le\MA_A(X)$. The inequality
$\MA_A(X)\le\Gam(X)$ follows from the definition of $\Gam(X)$.

Now take $A=H$ in \eqref{eq:energy-bottle}. Since
$T_H^X(x,y)=(x+y,x-y)$, the definition of $\JA_H(X)$ gives
\[
 \JA_H(X)=\sup_{x,y\in S_X}\min\{\|x+y\|,\|x-y\|\}=J(X).
\]
Dividing by $\|H\|_2=\sqrt2$ gives
$\widehat{\JA}_H(X)=J(X)/\sqrt2$. By
Proposition~\ref{prop:Hadamard}, $\MA_H(X)^2=\CNJ(X)$.
As $\MA_H(X)$ is nonnegative, $\MA_H(X)=\sqrt{\CNJ(X)}$.
Substituting these expressions into \eqref{eq:energy-bottle} gives
\[
 \frac{J(X)}{\sqrt2}\le\sqrt{\CNJ(X)}\le\Gam(X).
\]
Finally, Theorem~\ref{thm:universal} gives $\MA_A(X)\le\sqrt2$ for
every nonzero $A$. Taking the supremum over nonzero matrices yields
$\Gam(X)\le\sqrt2$, completing \eqref{eq:J-CNJ-Gamma}.
\end{proof}

Both sides of $J(X)/\sqrt2\le\sqrt{\CNJ(X)}$ in
\eqref{eq:J-CNJ-Gamma} are nonnegative. Squaring this inequality gives
\begin{equation}\label{eq:classical-J-CNJ}
 \frac{J(X)^2}{2}\le\CNJ(X).
\end{equation}
The inequality above also appears in
\cite[Theorem~3, Eq.~(3.1)]{KatoMaligrandaTakahashi2001}.
Stern proved $\CNJ(X)\le\CBM(X)\le2$ for $\dim X\ge2$
\cite[Theorem~3.4]{Stern2015}. The next proposition places $\Gam(X)^2$
between $\CNJ(X)$ and $\CBM(X)$.

\begin{proposition}\label{prop:old-constants}
For every real Banach space $X$,
\begin{equation}\label{eq:old-chain}
 \max\left\{\Cprime(X),\CZ(X),\frac{J(X)^2}{2}\right\}
 \le\CNJ(X)=\MA_H(X)^2
 \le\Gam(X)^2\le\CBM(X)\le2.
\end{equation}
\end{proposition}

\begin{proof}
Section~2 gives $\Cprime(X)\le\CNJ(X)$ and $\CZ(X)\le\CNJ(X)$.
Inequality \eqref{eq:classical-J-CNJ} gives the third bound.
The identity $\CNJ(X)=\MA_H(X)^2$ follows from
Proposition~\ref{prop:Hadamard}, and $\MA_H(X)\le\Gam(X)$ holds
by definition.

To prove $\Gam(X)^2\le\CBM(X)$, fix $A\ne0$ and a nonzero pair
$(x,y)\in X\oplus_2X$. Let $E=\spanop\{x,y\}$. The two output
coordinates of $T_A^X(x,y)$ are linear combinations of $x$ and $y$,
so $T_A^X$ restricts to $T_A^E$ on $E\oplus_2E$. Since $E$ has
the inherited norm, it follows that
\[
 \frac{\|T_A^X(x,y)\|_{2,X}}{\|A\|_2\|(x,y)\|_{2,X}}
 =\frac{\|T_A^E(x,y)\|_{2,E}}{\|A\|_2\|(x,y)\|_{2,E}}
 \le\MA_A(E).
\]
As $(x,y)\ne(0,0)$, the dimension of $E$ is either one or two.

\smallskip
\noindent\emph{Case 1. $\dim E=1$.}
The space $E$ is Hilbert, so $\MA_A(E)=1$ by the Hilbert space case
of Theorem~\ref{thm:universal}. If $\dim X=1$, then $\CBM(X)=1$
by convention. If $\dim X\ge2$, then $\CBM(X)\ge1$ by
\eqref{eq:CBM}, since the Banach--Mazur distance between
isomorphic normed spaces is at least one. In either case,
\[
 \MA_A(E)=1\le\sqrt{\CBM(X)}.
\]

\smallskip
\noindent\emph{Case 2. $\dim E=2$.}
Choose a linear isomorphism $S:E\to\ell_2^2$ and define
\[
 \widetilde S=S\oplus S:
 E\oplus_2E\longrightarrow\ell_2^2\oplus_2\ell_2^2,\qquad
 \widetilde S(u,v)=(Su,Sv).
\]
By linearity, $\widetilde S T_A^E=T_A^{\ell_2^2}\widetilde S$,
and consequently $T_A^E=\widetilde S^{-1}T_A^{\ell_2^2}\widetilde S$.
The direct sum norm gives $\|\widetilde S\|=\|S\|$ and
$\|\widetilde S^{-1}\|=\|S^{-1}\|$. For each equality, the upper
bound follows by applying the operator norm to both coordinates,
and the lower bound follows by restricting to pairs with second
coordinate zero. Since $\ell_2^2$ is Hilbert,
$\|T_A^{\ell_2^2}\|=\|A\|_2$. Hence
\[
 \|T_A^E\|
 \le\|\widetilde S^{-1}\|\,\|T_A^{\ell_2^2}\|\,\|\widetilde S\|
 =\|S^{-1}\|\,\|A\|_2\,\|S\|.
\]
Dividing by $\|A\|_2>0$ and taking the infimum over the
isomorphisms $S:E\to\ell_2^2$ yields
\[
 \MA_A(E)\le d_{\mathrm{BM}}(E,\ell_2^2)\le\sqrt{\CBM(X)}.
\]

Thus in both cases the normalized quotient for $(x,y)$ is at most
$\sqrt{\CBM(X)}$. Taking the supremum over nonzero pairs first
gives $\MA_A(X)\le\sqrt{\CBM(X)}$. Taking the supremum over
$A\ne0$ then gives $\Gam(X)\le\sqrt{\CBM(X)}$, and squaring
yields $\Gam(X)^2\le\CBM(X)$.

Finally, John's theorem \cite{John1948} gives
$d_{\mathrm{BM}}(E,\ell_2^2)\le\sqrt2$ for every two-dimensional
subspace $E$ of $X$. Squaring and taking the supremum in
\eqref{eq:CBM} proves $\CBM(X)\le2$ when $\dim X\ge2$.
For $\dim X=1$, the bound follows from $\CBM(X)=1$.
\end{proof}

\begin{remark}
On a nonzero Hilbert space, $\CNJ(X)=\Gam(X)^2=\CBM(X)=1$.
If $X$ contains an isometric copy of $\ell_1^2$ or $\ell_\infty^2$,
Theorem~\ref{thm:extreme-exact} gives $\MA_H(X)=\sqrt2$, and
\eqref{eq:old-chain} then gives $\CNJ(X)=\Gam(X)^2=\CBM(X)=2$.
The comparisons with $J(X)$ and $\Cprime(X)$ involve pairs of
unit vectors, so the equality cases need not coincide throughout
\eqref{eq:old-chain}.
\end{remark}

\section{Duality, subspaces and Hilbertian sums}

Passing to the dual space transposes the matrix in the formula for
$\MA_A$. Restricting to a subspace can only decrease $\MA_A$, and on a
Hilbertian sum its value is the supremum of the values on the summands.
We also bound how much $\MA_A$ can change under an isomorphism.

\begin{theorem}\label{thm:duality}
For every $A\ne0$,
\begin{equation}\label{eq:duality}
 \MA_A(X)=\MA_{A^T}(X^*).
\end{equation}
In particular,
\begin{equation}\label{eq:Gamma-dual}
 \Gam(X)=\Gam(X^*),\qquad
 \Lam(X;K)=\Lam(X^*;K).
\end{equation}
\end{theorem}

\begin{proof}
The dual of $X\oplus_2X$ is canonically $X^*\oplus_2X^*$, and
$(T_A^X)^*=T_{A^T}^{X^*}$.  Since an operator and its adjoint have the same
norm and $\|A^T\|_2=\|A\|_2$, \eqref{eq:duality} follows.  Transposition
preserves invertibility and the Euclidean condition number, so taking
suprema gives \eqref{eq:Gamma-dual}.
\end{proof}

\begin{proposition}\label{prop:subspace}
If $A\ne0$ and $Y$ is a nonzero closed subspace of $X$, then
\begin{equation}\label{eq:subspace}
 \MA_A(Y)\le\MA_A(X),\qquad
 \Gam(Y)\le\Gam(X).
\end{equation}
Moreover,
\begin{equation}\label{eq:finitely-determined}
 \MA_A(X)=\sup\{\MA_A(E):E\subseteq X,\ 1\le\dim E\le2\},
\end{equation}
and the analogous formula holds for $\Gam$.
\end{proposition}

\begin{proof}
The inherited norm identifies $Y\oplus_2Y$ isometrically with a
closed subspace of $X\oplus_2X$. Since both coordinates of
$T_A^X(u,v)$ are linear combinations of $u,v\in Y$, this subspace is
invariant and $T_A^Y$ is the restriction of $T_A^X$ to it. Hence
\[
 \|T_A^Y\|\le\|T_A^X\|.
\]
Dividing by $\|A\|_2>0$ gives $\MA_A(Y)\le\MA_A(X)$.
Taking the supremum over $A\ne0$ gives $\Gam(Y)\le\Gam(X)$.

For \eqref{eq:finitely-determined}, every subspace $E\subseteq X$ with
$1\le\dim E\le2$ is closed. Applying the preceding inequality to $E$
therefore gives
\[
 \sup_{\substack{E\subseteq X\\1\le\dim E\le2}}\MA_A(E)
 \le\MA_A(X).
\]
Conversely, take any nonzero pair $(x,y)\in X\oplus_2X$ and set
$E=\spanop\{x,y\}$. Then $1\le\dim E\le2$, and both output
coordinates of $T_A^X(x,y)$ lie in $E$. Since $E$ has the inherited
norm, the definition of the operator norm gives
\begin{align*}
 \frac{\|T_A^X(x,y)\|_{2,X}}{\|A\|_2\|(x,y)\|_{2,X}}
 &=\frac{\|T_A^E(x,y)\|_{2,E}}{\|A\|_2\|(x,y)\|_{2,E}}\\
 &\le \frac{\|T_A^E\|}{\|A\|_2}
 =\MA_A(E)
 \le\sup_{\substack{F\subseteq X\\1\le\dim F\le2}}\MA_A(F).
\end{align*}
Taking the supremum over all nonzero pairs $(x,y)$ yields the reverse
inequality. Combining the two inequalities proves
\eqref{eq:finitely-determined}.

The admissible subspaces $E$ in \eqref{eq:finitely-determined}
are the same for every nonzero $A$. Both iterated suprema below
therefore equal the supremum over all admissible pairs $(A,E)$.
\begin{align*}
 \Gam(X)
 &=\sup_{A\ne0}\MA_A(X)
 =\sup_{A\ne0}
   \sup_{\substack{E\subseteq X\\1\le\dim E\le2}}\MA_A(E)\\
 &=\sup_{\substack{E\subseteq X\\1\le\dim E\le2}}
   \sup_{A\ne0}\MA_A(E)
 =\sup_{\substack{E\subseteq X\\1\le\dim E\le2}}\Gam(E).
\end{align*}
\end{proof}

\begin{theorem}\label{thm:l2-sums}
Let $I$ be nonempty, let each $X_\alpha$ be a nonzero real Banach space, and put $X=(\bigoplus_{\alpha\in I}X_\alpha)_2$. For $A\ne0$ and $K\ge1$,
\begin{equation}\label{eq:l2-sum-M}
 \MA_A(X)=\sup_{\alpha\in I}\MA_A(X_\alpha)
\end{equation}
and consequently
\begin{equation}\label{eq:l2-sum-Gamma}
 \Gam(X)=\sup_{\alpha\in I}\Gam(X_\alpha),\qquad
 \Lam(X;K)=\sup_{\alpha\in I}\Lam(X_\alpha;K).
\end{equation}
\end{theorem}

\begin{proof}
Fix $A=(a_{ij})\ne0$. For $(x,y)\in X\oplus_2X$, write
$x=(x_\alpha)_{\alpha\in I}$ and $y=(y_\alpha)_{\alpha\in I}$.
The two components of $T_A^X(x,y)$ are
\[
 a_{11}x+a_{12}y=(a_{11}x_\alpha+a_{12}y_\alpha)_{\alpha\in I},
 \qquad
 a_{21}x+a_{22}y=(a_{21}x_\alpha+a_{22}y_\alpha)_{\alpha\in I}.
\]
Thus, grouping the two components with the same index $\alpha$,
the $\alpha$-coordinate pair of $T_A^X(x,y)$ is
$T_A^{X_\alpha}(x_\alpha,y_\alpha)$.

Set $M=\sup_{\alpha\in I}\|T_A^{X_\alpha}\|$. By
Theorem~\ref{thm:universal}, $M\le\|A^{\mathrm{abs}}\|_2$.
For arbitrary $I$, sums of nonnegative terms are understood
as suprema of finite subsums. The norm of the $\ell_2$-sum gives
\begin{align*}
 \|T_A^X(x,y)\|_{2,X}^2
 &=\sum_{\alpha\in I}
   \|T_A^{X_\alpha}(x_\alpha,y_\alpha)\|_{2,X_\alpha}^2\\
 &\le\sum_{\alpha\in I}
   \|T_A^{X_\alpha}\|^2\|(x_\alpha,y_\alpha)\|_{2,X_\alpha}^2\\
 &\le M^2\sum_{\alpha\in I}\|(x_\alpha,y_\alpha)\|_{2,X_\alpha}^2
 =M^2\|(x,y)\|_{2,X}^2.
\end{align*}
Taking square roots and then the supremum over unit pairs gives
$\|T_A^X\|\le M$.

For the reverse inequality, fix $\varepsilon>0$ and choose
$\alpha_0\in I$ such that $\|T_A^{X_{\alpha_0}}\|>M-\varepsilon$.
Since $M-\varepsilon$ is strictly below that operator norm,
there is a unit pair $(u,v)\in X_{\alpha_0}\oplus_2 X_{\alpha_0}$
with $\|T_A^{X_{\alpha_0}}(u,v)\|_{2,X_{\alpha_0}}>M-\varepsilon$.
Define $(x,y)\in X\oplus_2 X$ by
\[
 x_\alpha=\begin{cases}
 u,&\alpha=\alpha_0,\\
 0,&\alpha\ne\alpha_0,
 \end{cases}
 \qquad
 y_\alpha=\begin{cases}
 v,&\alpha=\alpha_0,\\
 0,&\alpha\ne\alpha_0.
 \end{cases}
\]
Then $\|(x,y)\|_{2,X}=\|(u,v)\|_{2,X_{\alpha_0}}=1$ and
\[
 \|T_A^X(x,y)\|_{2,X}
 =\|T_A^{X_{\alpha_0}}(u,v)\|_{2,X_{\alpha_0}}>M-\varepsilon.
\]
It follows that $\|T_A^X\|>M-\varepsilon$. Letting
$\varepsilon\downarrow0$ and combining both bounds gives
$\|T_A^X\|=\sup_{\alpha\in I}\|T_A^{X_\alpha}\|$.
Dividing by $\|A\|_2$ proves \eqref{eq:l2-sum-M}.

Taking the supremum in \eqref{eq:l2-sum-M} over nonzero matrices
and exchanging the order of the two suprema, we obtain
\begin{align*}
 \Gam(X)
 &=\sup_{A\ne0}\MA_A(X)
 =\sup_{A\ne0}\sup_{\alpha\in I}\MA_A(X_\alpha)\\
 &=\sup_{\alpha\in I}\sup_{A\ne0}\MA_A(X_\alpha)
 =\sup_{\alpha\in I}\Gam(X_\alpha).
\end{align*}
For $\Lam(X;K)$, take the supremum over invertible matrices
$A$ with $\cond(A)\le K$. The admissible set is the same for $X$
and each $X_\alpha$. Exchanging the suprema over this set and $I$
gives $\Lam(X;K)=\sup_{\alpha\in I}\Lam(X_\alpha;K)$.
\end{proof}

\begin{proposition}\label{prop:BM-stability}
If $S:X\to Y$ is an isomorphism, $A\ne0$, and
$d=\|S\|\|S^{-1}\|$, then
\begin{equation}\label{eq:BM-stability}
 d^{-1}\MA_A(Y)\le\MA_A(X)\le d\MA_A(Y).
\end{equation}
The same two-sided estimate holds for $\Gam$ and $\Lam(\,\cdot\,;K)$.
\end{proposition}

\begin{proof}
Define the induced isomorphism
\[
 \widetilde S=S\oplus S:X\oplus_2X\to Y\oplus_2Y,
 \qquad \widetilde S(x,y)=(Sx,Sy).
\]
For $(x,y)\in X\oplus_2X$,
\[
 \|\widetilde S(x,y)\|_{2,Y}^2
 =\|Sx\|_Y^2+\|Sy\|_Y^2
 \le\|S\|^2\|(x,y)\|_{2,X}^2.
\]
Thus $\|\widetilde S\|\le\|S\|$. The inverse is
$\widetilde S^{-1}=S^{-1}\oplus S^{-1}$, so the same estimate
gives $\|\widetilde S^{-1}\|\le\|S^{-1}\|$.
Linearity of $S$ gives, for every $(x,y)$,
\[
 \widetilde S T_A^X(x,y)
 =(a_{11}Sx+a_{12}Sy,a_{21}Sx+a_{22}Sy)
 =T_A^Y\widetilde S(x,y).
\]
Consequently $T_A^X=\widetilde S^{-1}T_A^Y\widetilde S$, and
\[
 \|T_A^X\|\le\|\widetilde S^{-1}\|\,\|T_A^Y\|\,
 \|\widetilde S\|\le d\|T_A^Y\|.
\]
Dividing by $\|A\|_2>0$ gives $\MA_A(X)\le d\MA_A(Y)$.
Applying the same argument to $S^{-1}:Y\to X$, whose distortion is
also $d$, gives $\MA_A(Y)\le d\MA_A(X)$. This proves
\eqref{eq:BM-stability}.

Taking suprema over $A\ne0$ in both inequalities gives
\[
 d^{-1}\Gam(Y)\le\Gam(X)\le d\Gam(Y).
\]
For $K\ge1$, taking suprema over the invertible matrices satisfying
$\cond(A)\le K$ gives
\[
 d^{-1}\Lam(Y;K)\le\Lam(X;K)\le d\Lam(Y;K).
\]
\end{proof}

\section{Exact calculations in classical spaces}

We use complex interpolation to bound $\MA_A$ on $L_p$ and $\ell_p$
for each matrix $A$. To calculate the global profile, we use the
Hadamard matrix and construct pairs from disjointly supported vectors
to attain the upper bound.

\begin{theorem}\label{thm:Lp}
Let $1\le p\le\infty$, let $A\ne0$, and assume the $L_p$ space under
consideration is nonzero.  Then
\begin{equation}\label{eq:Lp-individual}
 \MA_A(L_p)\le
 \left(\frac{\|A^{\mathrm{abs}}\|_2}{\|A\|_2}\right)^{|2/p-1|},
\end{equation}
where $1/\infty=0$.  The same estimate holds for $\ell_p$.
If the space contains two disjointly supported nonzero vectors, then
\begin{equation}\label{eq:Lp-profile}
 \Gam(L_p)=\Lam(L_p;K)=2^{|1/p-1/2|}
 \qquad(K\ge1),
\end{equation}
and likewise for $\ell_p$ of dimension at least two.
\end{theorem}

\begin{proof}
\leavevmode\par\nopagebreak[4]
\noindent\textbf{The upper estimate \eqref{eq:Lp-individual}.}
\par\nopagebreak[4]
We first work over the complex scalars and set
$E_p=L_p(\mu;\mathbb C)\oplus_2 L_p(\mu;\mathbb C)$.
The Hilbert space estimate on $E_2$ gives
$\|T_A\|_{E_2\to E_2}\le\|A\|_2$. The triangle-inequality
argument of Theorem~\ref{thm:universal} also applies to complex
spaces, giving
\[
 \|T_A\|_{E_1\to E_1}\le\|A^{\mathrm{abs}}\|_2,\qquad
 \|T_A\|_{E_\infty\to E_\infty}\le\|A^{\mathrm{abs}}\|_2.
\]

First assume $\mu(\Omega)<\infty$, so simple functions of finite-measure
support are dense in $L_\infty(\mu;\mathbb C)$.
For either interpolation couple, the same matrix formula defines a
complex-linear operator on the sum of the endpoint spaces.

For $1<p<2$, choose $\theta\in(0,1)$ so that
\[
 \frac1p=\frac{1-\theta}{1}+\frac{\theta}{2}
 =1-\frac{\theta}{2},\qquad \theta=2-\frac2p.
\]
Scalar $L_p$ interpolation followed by Bochner $L_2$ interpolation
on the two-point counting measure space gives
$E_p=[E_1,E_2]_\theta$ isometrically
\cite[Theorems~5.1.1 and~5.1.2]{BerghLofstrom1976}.
Applying the operator interpolation theorem
\cite[Theorem~4.1.2]{BerghLofstrom1976} to $T_A$ gives
\begin{align*}
 \|T_A\|_{E_p\to E_p}
 &\le \|T_A\|_{E_1\to E_1}^{1-\theta}
       \|T_A\|_{E_2\to E_2}^{\theta}\\
 &\le \|A^{\mathrm{abs}}\|_2^{1-\theta}\|A\|_2^\theta.
\end{align*}
Here $1-\theta=2/p-1=|2/p-1|$ and
$\theta=1-|2/p-1|$.

For $2<p<\infty$, the interpolation couple is $(E_2,E_\infty)$.
The same interpolation identities give $E_p=[E_2,E_\infty]_\theta$
isometrically, where
\[
 \frac1p=\frac{1-\theta}{2}+\theta\cdot 0
 =\frac{1-\theta}{2},\qquad \theta=1-\frac2p.
\]
Applying the operator theorem to this couple gives
\begin{align*}
 \|T_A\|_{E_p\to E_p}
 &\le \|T_A\|_{E_2\to E_2}^{1-\theta}
       \|T_A\|_{E_\infty\to E_\infty}^{\theta}\\
 &\le \|A\|_2^{1-\theta}\|A^{\mathrm{abs}}\|_2^\theta.
\end{align*}
Here $\theta=1-2/p=|2/p-1|$ and $1-\theta=1-|2/p-1|$.

Combining the two cases with the direct endpoint bounds at
$p=1,2,\infty$ yields
\[
 \|T_A\|_{E_p\to E_p}
 \le \|A\|_2^{1-|2/p-1|}\,\|A^{\mathrm{abs}}\|_2^{|2/p-1|}.
\]
For $1<p<\infty$, simple functions supported on sets of finite measure
are dense in $L_p(\mu;\mathbb C)$. For any pair of such functions,
both coordinates of its image under $T_A$ vanish outside the union
of their supports. The finite-measure estimate thus applies to the
pair with the same constant. By density and continuity of $T_A$,
the estimate extends to $E_p$ for arbitrary $\mu$. The cases
$p=1,2,\infty$ follow from the direct endpoint estimates above.

The pairs of real-valued functions form an isometric subspace of
$E_p$ that is invariant under $T_A$, because the entries of $A$
are real. The operator norm of this restriction is bounded by
the norm on $E_p$. Dividing the estimate by $\|A\|_2>0$ therefore
proves \eqref{eq:Lp-individual} for the real $L_p$ space. Counting
measure gives the same estimate for $\ell_p$.

\smallskip
\noindent\textbf{The exact value \eqref{eq:Lp-profile}.}
\par\nopagebreak[4]
Since $\|A^{\mathrm{abs}}\|_2/\|A\|_2\le\sqrt2$ by
Theorem~\ref{thm:universal}, the upper estimate gives
\[
 \MA_A(L_p)\le(\sqrt2)^{|2/p-1|}
 =2^{|1/p-1/2|}
 \qquad (A\ne0).
\]
Taking the supremum over $A\ne0$ gives
$\Gam(L_p)\le 2^{|1/p-1/2|}$.

For the reverse bound, take disjointly supported nonzero real
functions in $L_p$ and normalize them to unit vectors $u,v$.
We use $T_H^{L_p}(x,y)=(x+y,x-y)$ and $\|H\|_2=\sqrt2$.
For any unit vectors $x,y\in L_p$, the norm of the input pair is
$\sqrt2$, so the normalized quotient gives
\[
 \MA_H(L_p)\ge
 \frac{\|T_H^{L_p}(x,y)\|_{2,L_p}}
 {\|H\|_2\|(x,y)\|_{2,L_p}}
 =\frac{(\|x+y\|_p^2+\|x-y\|_p^2)^{1/2}}{2}.
\]

\emph{Case $1\le p\le2$.}
Take $x=u$ and $y=v$. Disjointness gives
$|u\pm v|^p=|u|^p+|v|^p$ almost everywhere, hence
\[
 \|u+v\|_p^p=\|u-v\|_p^p=\|u\|_p^p+\|v\|_p^p=2.
\]
Both output norms are $2^{1/p}$. Thus
\[
 \MA_H(L_p)\ge
 \frac{(2\cdot 2^{2/p})^{1/2}}{2}
 =2^{1/p-1/2}=2^{|1/p-1/2|}.
\]

\emph{Case $2<p\le\infty$.}
Set $x=2^{-1/p}(u+v)$ and $y=2^{-1/p}(u-v)$.
If $p<\infty$, disjointness gives
$\|u\pm v\|_p=2^{1/p}$, so both $x$ and $y$ are unit vectors.
If $p=\infty$, the disjoint supports give
$\|u\pm v\|_\infty=\max\{\|u\|_\infty,\|v\|_\infty\}=1$,
and the factor $2^{-1/p}$ equals one. Thus $x,y$ are unit
vectors at this endpoint as well. Moreover,
\[
 x+y=2^{1-1/p}u,\qquad x-y=2^{1-1/p}v.
\]
Both output norms are $2^{1-1/p}$, so
\[
 \MA_H(L_p)\ge
 \frac{(2\cdot 2^{2(1-1/p)})^{1/2}}{2}
 =2^{1/2-1/p}=2^{|1/p-1/2|}.
\]

By Proposition~\ref{prop:Hadamard} and the calculation of
$\cond(H)$ in its proof, $H$ is admissible in the definition of
$\Lam(L_p;K)$ for every $K\ge1$. We therefore have
\[
 2^{|1/p-1/2|}\le\MA_H(L_p)
 \le\Lam(L_p;K)\le\Gam(L_p)\le 2^{|1/p-1/2|},
\]
which proves \eqref{eq:Lp-profile}. For $\ell_p$ of dimension at
least two, two different coordinate unit vectors give the same
disjoint-support computations and the same equalities.
\end{proof}

The vector pairs that attain the Hadamard amplification bound also
attain the bound for the James constant. This gives both classical
constants in terms of the global profile.

\begin{corollary}\label{cor:Lp-J}
For $1\le p\le\infty$ and $q$ conjugate to $p$,
suppose that $L_p$ contains two disjointly supported nonzero
functions. Then
\begin{align}
 J(L_p)&=\max\{2^{1/p},2^{1/q}\},\label{eq:Lp-J}\\
 \CNJ(L_p)&=2^{|2/p-1|},\label{eq:Lp-CNJ}\\
 J(L_p)&=\sqrt2\,\Gam(L_p),\qquad
 \CNJ(L_p)=\Gam(L_p)^2.\label{eq:Lp-all}
\end{align}
The same formulas hold for $\ell_p$ of dimension at least two.
\end{corollary}

\begin{proof}
The Hadamard test pairs in the proof of Theorem~\ref{thm:Lp} give
$\MA_H(L_p)\ge 2^{|1/p-1/2|}$. Since
$\MA_H(L_p)\le\Gam(L_p)$ and
$\Gam(L_p)=2^{|1/p-1/2|}$ by \eqref{eq:Lp-profile}, equality holds.
Proposition~\ref{prop:Hadamard} therefore gives
\[
 \CNJ(L_p)=\MA_H(L_p)^2=\Gam(L_p)^2
 =2^{2|1/p-1/2|}=2^{|2/p-1|},
\]
which proves \eqref{eq:Lp-CNJ} and the second identity in
\eqref{eq:Lp-all}. For the classical computation of this constant,
see \cite{Clarkson1937,KatoTakahashi1998}.

To obtain \eqref{eq:Lp-J} and the first identity in \eqref{eq:Lp-all},
we use the upper bound $J(L_p)\le\sqrt2\,\Gam(L_p)$ from
\eqref{eq:J-CNJ-Gamma} and the unit pairs constructed in
Theorem~\ref{thm:Lp}.

\emph{Case $1\le p\le2$.}
Here $1/q\le 1/2\le 1/p$, so
\[
 \sqrt2\,\Gam(L_p)=2^{1/2}2^{1/p-1/2}
 =2^{1/p}=\max\{2^{1/p},2^{1/q}\}.
\]
The disjointly supported unit pair $(u,v)$ satisfies
$\|u+v\|_p=\|u-v\|_p=2^{1/p}$. Thus $J(L_p)\ge 2^{1/p}$,
and the upper bound gives equality.

\emph{Case $2<p\le\infty$.}
Here $1/p\le 1/2\le 1/q$, so
\[
 \sqrt2\,\Gam(L_p)=2^{1/2}2^{1/2-1/p}
 =2^{1-1/p}=2^{1/q}=\max\{2^{1/p},2^{1/q}\}.
\]
For the unit pair $x=2^{-1/p}(u+v)$ and $y=2^{-1/p}(u-v)$,
both $\|x+y\|_p$ and $\|x-y\|_p$ equal $2^{1/q}$.
Thus $J(L_p)\ge 2^{1/q}$, and the upper bound again gives equality.
This proves \eqref{eq:Lp-J} and $J(L_p)=\sqrt2\,\Gam(L_p)$.
For the classical values of the James constant on $\ell_p$
and $L_p[0,1]$, see Gao and Lau \cite[Theorems~3.1 and~3.2]{GaoLau1990}.
The formula for $L_p(\mu)$ of dimension at least two is recorded in
\cite[p.~280]{KatoMaligrandaTakahashi2001}.

For $\ell_p$ of dimension at least two, the same calculations apply
with two distinct coordinate unit vectors as the disjointly
supported vectors. The endpoint cases use $1/\infty=0$.
\end{proof}

Table~\ref{tab:values} lists $J(X)$, $\CNJ(X)$ and $\Gam(X)$ for
the spaces considered above. In every $L_p$ row the space is assumed to contain two
disjointly supported nonzero vectors. In the $C(K)$ row, $K$ is a compact
Hausdorff space with at least two points.

\begin{table}[htbp]
\centering
\caption{Values of the classical constants and the global profile.}
\label{tab:values}
\renewcommand{\arraystretch}{1.18}
\begin{tabular}{@{}>{\raggedright\arraybackslash}p{0.28\textwidth}
                >{\centering\arraybackslash}p{0.18\textwidth}
                >{\centering\arraybackslash}p{0.18\textwidth}
                >{\centering\arraybackslash}p{0.18\textwidth}@{}}
\toprule
Space $X$ & $J(X)$ & $\CNJ(X)$ & $\Gam(X)$\\
\midrule
Hilbert, $\dim X\ge2$ & $\sqrt2$ & $1$ & $1$\\
$\ell_p,L_p$, $1<p<\infty$ &
$2^{\max\{1/p,1-1/p\}}$ & $2^{|2/p-1|}$ &
$2^{|1/p-1/2|}$\\
$\ell_1,L_1$ & $2$ & $2$ & $\sqrt2$\\
$c_0,\ell_\infty,L_\infty,C(K)$ & $2$ & $2$ & $\sqrt2$\\
\bottomrule
\end{tabular}
\end{table}

\begin{example}
Let $X=(\bigoplus_{n\ge2}\ell_{p_n}^n)_2$, where $p_n\in[1,\infty]$.
Theorem~\ref{thm:l2-sums} gives the profile of the direct sum, and
Theorem~\ref{thm:Lp} gives the profile of each summand. Thus
\[
 \Gam(X)=\sup_{n\ge2}\Gam(\ell_{p_n}^n)
 =\sup_{n\ge2}2^{|1/p_n-1/2|}.
\]
By Theorem~\ref{thm:global-profile}, $X$ is uniformly nonsquare exactly
when this supremum is strictly below $\sqrt2$.

If some $p_n$ equals $1$ or $\infty$, its summand has profile $\sqrt2$.
The supremum is also $\sqrt2$ if a subsequence of $(p_n)$ tends to an
endpoint. These cases cover failure of the common bounds below. In fact, $p_{n_k}\to1$ gives $1/p_{n_k}\to1$, and
$p_{n_k}\to\infty$ gives $1/p_{n_k}\to0$. In both cases
$2^{|1/p_{n_k}-1/2|}\to\sqrt2$.

If $1<\inf_{n\ge2}p_n\le\sup_{n\ge2}p_n<\infty$, there is a single
$\delta>0$ such that $\delta\le1/p_n\le1-\delta$ for all $n\ge2$.
Consequently,
\[
 \Gam(\ell_{p_n}^n)\le 2^{1/2-\delta}<\sqrt2
 \qquad(n\ge2).
\]
Taking the supremum preserves this common upper bound, so $X$ is
uniformly nonsquare.

Thus $X$ is uniformly nonsquare if and only if the exponents $p_n$
are uniformly bounded away from $1$ and $\infty$, that is,
$1<\inf_{n\ge2}p_n\le\sup_{n\ge2}p_n<\infty$.

The choice $p_n=n$ illustrates the case $p_n\to\infty$. Each summand is
uniformly nonsquare, and their profiles converge to $\sqrt2$. The direct
sum then has profile $\sqrt2$.

\end{example}

\section{Convexity, smoothness and fixed points}

The modulus of convexity gives bounds on the James constant, and
the modulus of smoothness at $1$ is bounded in terms of
$\CNJ(X)$. We then use $\Gam(X)$ to give sufficient conditions for the
fixed point property and for uniform normal structure.

\begin{proposition}\label{prop:delta}
For every $\varepsilon\in(0,2)$,
\begin{equation}\label{eq:J-delta-bound}
 J(X)\le
 \max\{\varepsilon,2(1-\delta_X(\varepsilon))\}.
\end{equation}
Consequently,
\begin{equation}\label{eq:delta-operator}
 \frac{1}{\sqrt2}
 \sup_{0<\varepsilon<2}
 \min\{\varepsilon,2(1-\delta_X(\varepsilon))\}
 \le\MA_H(X)\le\Gam(X).
\end{equation}
More precisely, $J(X)$ is given by \eqref{eq:J-delta}.
\end{proposition}

\begin{proof}
To prove \eqref{eq:J-delta-bound}, fix $\varepsilon\in(0,2)$
and take $x,y\in S_X$. If $\|x-y\|<\varepsilon$, then
\[
 \min\{\|x+y\|,\|x-y\|\}\le\|x-y\|<\varepsilon.
\]
If $\|x-y\|\ge\varepsilon$, the definition \eqref{eq:delta} gives
$\delta_X(\varepsilon)\le1-\|(x+y)/2\|$. It follows that
\[
 \min\{\|x+y\|,\|x-y\|\}
 \le\|x+y\|\le2\bigl(1-\delta_X(\varepsilon)\bigr).
\]
Thus the minimum is at most
$\max\{\varepsilon,2(1-\delta_X(\varepsilon))\}$ for every
$x,y\in S_X$. Taking the supremum over $x,y\in S_X$ proves
\eqref{eq:J-delta-bound}.

For the lower estimate in \eqref{eq:delta-operator}, fix
$\varepsilon\in(0,2)$ and $\eta>0$. The pair $(u,-u)$ for any $u\in S_X$ has distance $2$, so the
defining set in \eqref{eq:delta} is nonempty. By the definition
of the infimum,
there exist $x,y\in S_X$ such that $\|x-y\|\ge\varepsilon$ and
\[
 1-\left\|\frac{x+y}{2}\right\|<\delta_X(\varepsilon)+\eta.
\]
The latter inequality gives $\|x+y\|>2(1-\delta_X(\varepsilon))-2\eta$.
By the definition of $J(X)$,
\[
 J(X)\ge\min\{\|x+y\|,\|x-y\|\}
 \ge\min\{\varepsilon,2(1-\delta_X(\varepsilon))-2\eta\}.
\]
Letting $\eta\downarrow0$ and using continuity of the minimum gives
\[
 J(X)\ge\min\{\varepsilon,2(1-\delta_X(\varepsilon))\}.
\]
Taking the supremum over $\varepsilon\in(0,2)$, dividing by $\sqrt2$,
and applying \eqref{eq:J-CNJ-Gamma} and \eqref{eq:Hadamard-CNJ}
give \eqref{eq:delta-operator}.

To verify \eqref{eq:J-delta}, observe that the preceding bound
gives $J(X)\ge\varepsilon$ whenever
$\delta_X(\varepsilon)\le1-\varepsilon/2$. Conversely, if
$0<\varepsilon<J(X)$, there are $x,y\in S_X$ with
$\min\{\|x+y\|,\|x-y\|\}>\varepsilon$. This pair is admissible in
\eqref{eq:delta}, so
$\delta_X(\varepsilon)\le1-\|x+y\|/2<1-\varepsilon/2$.
Since $0\le J(X)\le2$, these two implications yield
\eqref{eq:J-delta}, including the case $J(X)=0$.
\end{proof}

\begin{proposition}\label{prop:rho}
For every real Banach space $X$,
\begin{equation}\label{eq:rho-bound}
 1+\rho_X(1)\le\sqrt{2\CNJ(X)}
 =\sqrt2\,\MA_H(X)\le\sqrt2\,\Gam(X).
\end{equation}
\end{proposition}

\begin{proof}
Take $x,y\in S_X$ and put $a=\|x+y\|$ and $b=\|x-y\|$.
By the Cauchy--Schwarz inequality applied to $(a,b)$ and $(1,1)$,
\[
 \frac{a+b}{2}\le\left(\frac{a^2+b^2}{2}\right)^{1/2}.
\]
Since $\|x\|=\|y\|=1$, the defining quotient in
\eqref{eq:CNJ-intro} gives
\[
 \frac{a^2+b^2}{4}
 =\frac{\|x+y\|^2+\|x-y\|^2}{2(\|x\|^2+\|y\|^2)}
 \le\CNJ(X).
\]
Consequently $(a+b)/2\le\sqrt{2\CNJ(X)}$. Taking the supremum over
$x,y\in S_X$ gives $1+\rho_X(1)\le\sqrt{2\CNJ(X)}$ by \eqref{eq:rho}.
Proposition~\ref{prop:Hadamard} gives
$\sqrt{\CNJ(X)}=\MA_H(X)\le\Gam(X)$, which completes
\eqref{eq:rho-bound}.
\end{proof}

By Theorem~\ref{thm:global-profile}, a bound $\Gam(X)<\sqrt2$ gives
uniform nonsquareness. The following theorem uses this to obtain the
fixed point property for nonexpansive mappings.

\begin{theorem}\label{thm:FPP}
If
\begin{equation}\label{eq:FPP-threshold}
 \Gam(X)<\sqrt2,
\end{equation}
then $X$ has the fixed point property for nonexpansive mappings: every
nonexpansive self-map of a nonempty bounded closed convex subset of $X$ has
a fixed point.
\end{theorem}

\begin{proof}
By Theorem~\ref{thm:global-profile}, condition
\eqref{eq:FPP-threshold} is equivalent to uniform nonsquareness.  The fixed
point property of uniformly nonsquare Banach spaces was proved by
Garc\'ia-Falset, Llorens-Fuster and Mazcu\~nan-Navarro
\cite{GarciaFalset2006}.  Applying their theorem gives the conclusion.
\end{proof}

\begin{corollary}\label{cor:UNS}
If
\begin{equation}\label{eq:uniform-normal}
 \Gam(X)<\frac{\sqrt5}{2},
\end{equation}
then $X$ has uniform normal structure and hence the fixed point property.
\end{corollary}

\begin{proof}
If $\dim X=1$, the midpoint of any bounded nondegenerate interval in
$X$ is at distance at most half its diameter from every point of the
interval. Hence $X$ has uniform normal structure. Every nonempty bounded
closed convex subset of $X$ is a closed interval, on which a nonexpansive
self-map has a fixed point by the intermediate value theorem.

Now suppose $\dim X\ge2$. By Proposition~\ref{prop:Hadamard},
$\CNJ(X)=\MA_H(X)^2\le\Gam(X)^2<5/4$.  The theorem of Kato, Maligranda
and Takahashi \cite[Corollary~3]{KatoMaligrandaTakahashi2001} implies uniform normal
structure under the condition $\CNJ(X)<5/4$.
Since $\CNJ(X)<2$, \eqref{eq:uns-equiv} gives uniform nonsquareness,
and James's theorem \cite{James1964} gives reflexivity. Every nonempty
bounded closed convex subset of $X$ is therefore weakly compact.
It has normal structure because $X$ has uniform normal structure,
so Kirk's theorem \cite{Kirk1965} gives a fixed point for every
nonexpansive self-map of such a subset.
\end{proof}

\begin{remark}
The condition $\Gam(X)<\sqrt2$ gives the fixed point property via
uniform nonsquareness. The smaller threshold $\sqrt5/2$ in
Corollary~\ref{cor:UNS} also guarantees uniform normal structure.
\end{remark}

\section{Further structural consequences}

For Hilbertian sums, uniform nonsquareness requires a common
upper bound strictly below $\sqrt2$ for the profiles of the summands. A uniform bound on the
Banach--Mazur distances of the two-dimensional subspaces to $\ell_2^2$
also gives a bound on $\Gam(X)$.

\begin{proposition}\label{prop:family-gap}
Let $I$ be nonempty and $(X_\alpha)_{\alpha\in I}$ be a family of nonzero real Banach spaces.  The Hilbertian
sum $X=(\bigoplus_\alpha X_\alpha)_2$ is uniformly nonsquare if and only if
there is $\eta>0$ such that
\begin{equation}\label{eq:family-gap}
 \Gam(X_\alpha)\le\sqrt2-\eta
 \qquad(\alpha\in I).
\end{equation}
\end{proposition}

\begin{proof}
Suppose first that $X$ is uniformly nonsquare.
Theorem~\ref{thm:global-profile} gives $\Gam(X)<\sqrt2$, and
Theorem~\ref{thm:l2-sums} gives
$\Gam(X)=\sup_{\alpha\in I}\Gam(X_\alpha)$.
Thus the common choice $\eta=\sqrt2-\Gam(X)>0$ satisfies
\[
 \Gam(X_\alpha)\le\Gam(X)=\sqrt2-\eta
 \qquad(\alpha\in I).
\]

Conversely, suppose that \eqref{eq:family-gap} holds for some
$\eta>0$. Taking the supremum over $\alpha\in I$ and applying
Theorem~\ref{thm:l2-sums}, we obtain
\[
 \Gam(X)=\sup_{\alpha\in I}\Gam(X_\alpha)
 \le\sqrt2-\eta<\sqrt2.
\]
Theorem~\ref{thm:global-profile} now implies that $X$ is uniformly
nonsquare.
\end{proof}

\begin{proposition}
\label{prop:near-Hilbert}
Let $\varepsilon\ge0$. If $\Gam(X)\le1+\varepsilon$, then every two-dimensional subspace $E$ of
$X$ satisfies
\begin{equation}\label{eq:near-Hilbert-oneway}
 \MA_A(E)\le1+\varepsilon
 \qquad(A\ne0),
\end{equation}
and
\begin{equation}\label{eq:near-Hilbert-classical}
 \CNJ(E)\le(1+\varepsilon)^2,\qquad
 J(E)\le\sqrt2(1+\varepsilon).
\end{equation}
Conversely, a uniform Banach--Mazur estimate
$d_{\mathrm{BM}}(E,\ell_2^2)\le1+\varepsilon$ for all two-dimensional
$E\subset X$ implies $\Gam(X)\le1+\varepsilon$.
\end{proposition}

\begin{proof}
The first assertion follows from subspace monotonicity.  The classical
bounds follow from \eqref{eq:J-CNJ-Gamma}.
For the converse, the Banach--Mazur hypothesis gives
$\CBM(X)\le(1+\varepsilon)^2$ by \eqref{eq:CBM}. The same bound
holds in dimension one by the convention $\CBM(X)=1$.
Proposition~\ref{prop:old-constants} gives
$\Gam(X)^2\le\CBM(X)\le(1+\varepsilon)^2$, and hence
$\Gam(X)\le1+\varepsilon$.
\end{proof}

\begin{remark}
The Banach--Mazur estimate in the hypothesis bounds $\Gam(X)$
from above through $\Gam(X)^2\le\CBM(X)$. To bound the
Banach--Mazur distances from a bound on $\Gam(X)$, an estimate in
the opposite direction would be needed.
\end{remark}

\section*{Acknowledgements}
\addcontentsline{toc}{section}{Acknowledgements}

Research supported by the National Natural Science Foundation of China (Nos.12471130; 11771126).

\end{document}